# THE OBSERVABILITY CONCEPT IN A CLASS OF HYBRID CONTROL SYSTEMS

ECATERINA OLTEAN

In the discrete modeling approach for hybrid control systems, the continuous plant is reduced to a discrete event approximation, called the DES-plant, that is governed by a discrete event system, representing the controller. The observability of the DES-plant model is crucial for the synthesis of the controller and for the proper closed loop evolution of the hybrid control system. Based on a version of the framework for hybrid control systems proposed by Antsaklis, the paper analysis the relation between the properties of the cellular space of the continuous plant and a mechanism of plant-symbols generation, on one side, and the observability of the DES-plant automaton on the other side. Finally an observable discrete event abstraction of the continuous double integrator is presented.

## 1. INTRODUCTION AND MOTIVATION

In a class of hybrid control systems (HCS), a discrete event system (DES), representing the controller, interacts, through an interface, with a purely continuous plant (fig.1). This framework, sometimes called the supervisor approach for HCS, was first proposed by Antsaklis and his co-workers [1], [2], and also by Nerode and Kohn [5]. The controller is modeled as a Moore machine. The plant coupled to the interface behaves like a DES, called the DES-plant, represented as a nondeterministic automaton. The plant's continuous state space is partitioned, by a finite set of smooth hypersurfaces, into open cells; each cell corresponds to a discrete state of the DES–plant automaton. A plant-event occurs, and a corresponding plant-symbol is sent through the interface, whenever a hypersurface is crossed, in a certain direction, by the continuous state trajectory. A plant-symbol labels the open halfspace in which the plant's state just enters. The interface also converts sequences of control-symbols into an input signal for the plant.

The only information received by the controller is the sequence of plant-symbols, generated when the continuous state trajectory crosses the hypersurfaces and evolves, from cell to cell, in the partitioned state space. On the other hand, the control objective is usually represented as a desired string of open cells, corresponding to a discrete evolution of the DES-plant; the continuous trajectory has to reside within the desired string of cells. The observability of DES-plant's evolution is crucial for the proper closed loop behaviour of the HCS. Intuitively, the observability property must ensure that an initial discrete state $p_0$ of the DES-plant model and an admissible sequence of plant-symbols correspond, together, to a unique discrete evolution of the DES-plant automaton, starting in $p_0$.

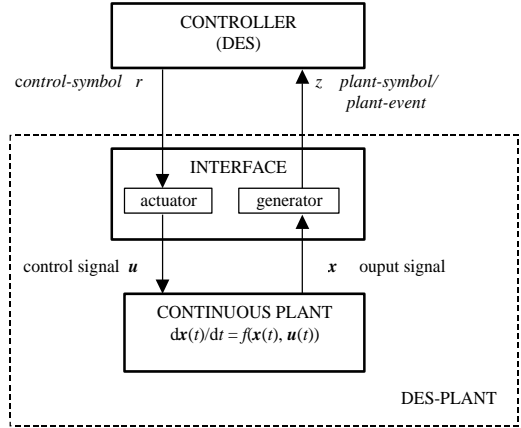

Fig. 1 - The architecture of a hybrid control system.

Based on a version of the Antsaklis framework, studied in [4], the paper presents the relation between the properties of the adjacent cells, the mechanism of plant-symbols generation and the observability of the DES-plant. First some basic definitions and properties concerning the cellular space are presented in section 2 and a mechanism for plant-symbols generation is proposed in section 3. Then the observability concept is introduced, in section 4, for a DES-plant model which is slightly different from the one first proposed in the Antsaklis framework. It is shown that a property of the cellular space ensures, under certain assumptions, the observability of the DES-plant automaton. Finally, a simple example of an observable DES-plant automaton, associated to the continuous double integrator, is discussed.

## 2. THE CONTINUOUS STATE SPACE PARTITION AND THE CELLULAR SPACE

### 2.1 THE CONTINUOUS PLANT AND THE ACTUATOR

Define $\mathbb{R}^+ = [0, \infty)$ the continuous time domain, $t \in \mathbb{R}^+$ the continuous time variable, $\mathbb{N} = \{0, 1, \ldots\}$ the logical time domain and $k \in \mathbb{N}$ the logical time variable. The *continuous plant* is modeled as a continuous time invariant dynamic system

$$(\text{ODE}) \quad \dot{x}(t) = f(x(t), u(t)) \tag{1}$$

where $x(t) \in X \subseteq \mathbb{R}^n$ and $\mathbf{u}(t) \in U \subset \mathbb{R}^m$ are the state and control vectors respectively, at the moment $t \in \mathbb{R}^+$. $X$ is the *continuous state space*. $U = \{\mathbf{u}_1, \ldots, \mathbf{u}_M\}$ is the discrete set of admissible control values and $R = \{r_1, \ldots, r_M\}$ is the *alphabet of control-symbols*. Consider the bijective map $act : R \to U$, $act(r_m) = \mathbf{u}_m$, $\forall \, m \in I_M = \{1, \ldots, M\}$. The *actuator* converts a sequence of control-symbols $w_r = r(0), r(1), \ldots, r(k), \ldots, r(k) \in R, \forall \, k \in \mathbb{N}$ into a piecewise constant control signal

$$\mathbf{u}(t) = \sum_{k \geq 0} act(r(k)) \cdot I(t, t_c(k), t_c(k+1)) \qquad (2)$$

where $t_c(k) \in \mathbb{R}^+$ is the moment when $r(k)$ is received from the DES controller, $t_c(k) < t_c(k+1), \forall \, k \in \mathbb{N}$, and $I : \mathbb{R}^+ \times \mathbb{R}^+ \times \mathbb{R}^+ \to \{0, 1\}$, $I(t, t_1, t_2) = 1$, if $t_1 \leq t < t_2$ and $I(t, t_1, t_2) = 0$, else. When a control-symbol $r_m \in R$ is received, the control signal $\mathbf{u}(.)$ instantly switches to the value $act(r_m) \in U$ and remains there till a new control-symbols occurs.

## 2.2 THE CELLULAR SPACE

Consider $N \geq 1$ a natural number and $I_N = \{1, \ldots, N\}$ an index set.

*Definition 1*. A) A *partition of the continuous state space* $X \subseteq \mathbb{R}^n$ is a set of $N$ indexed functionals, $S_h^N = \{h_i \mid h_i : X \to \mathbb{R}, i \in I_N\}$, with the following properties : $\forall \, h_i \in S_h^N$ a1) $h_i \in C^1$ in $X$ and a2) $\text{Ker}(h_i) = \{x \in X \mid h_i(x) = 0\}$ is a nonsingular $(n-1)$-dimensional hypersurface ( i.e. $\text{grad}_x(h_i) \neq 0, \forall \, x \in \text{Ker}(h_i)$ ).

B) $\forall \, h_i \in S_h^N$, $h_i$ separates $X$ into two open halfspaces $H_i^- = \{x \in X \mid h_i(x) < 0\}$ and $H_i^+ = \{x \in X \mid h_i(x) > 0\}$. $S_h^N$ defines the collection of $2N$ halfspaces $SH = \{H_1^-, H_1^+, \ldots, H_N^-, H_N^+\}$, called *event-topology*. □

In the Antsaklis framework, the following equivalence relation permits the building of the DES-plant's discrete state space.

*Definition 2*. [1] Consider $X \subseteq \mathbb{R}^n$, $DX = X \setminus Fr$, with $Fr = \bigcup_{i=1}^{N} \text{Ker}(h_i)$. The *equivalence relation* induced by $S_h^N$ is $Rel \subset DX \times DX$, defined by :

$$(x_1, x_2) \in Rel \Leftrightarrow h_i(x_1) h_i(x_2) > 0, \forall \, h_i \in S_h^N \qquad (3)$$

□

*Definition 3*. A) The *cellular space* $DX_{/Rel} = C = \{c_1, \ldots, c_Q\}$ is the set of all classes of equivalence of the relation $Rel$, where $c_q \in C$ is a cell, $I_Q = \{1, \ldots, Q\}$ is the index set of $C$, with $Q = \text{card}(C)$.

B) The *alphabet of discrete states* (or the discrete state space of the DES-plant) is a set of $Q$ distinct indexed symbols $P = \{p_1, \ldots, p_Q\}$. The map $et : C \to P$, $et(c_q) = p_q, \forall \, q \in I_Q$ is the *label-function* of $C$.

C) The *canonical sujection* with respect to $C$ is the map $sur_C : DX \to C$ defined s.t. $\forall \, x \in DX$, $sur_C(x) = c_q \Leftrightarrow x \in c_q$. □

Consider the function *sgn* : $\mathbb{R} \to \{-1, 0, 1\}$, $sgn(y) = -1$, if $y < 0$, $sgn(y) = 0$, if $y = 0$ and $sgn(y) = 1$, if $y > 0$. The next definition, introduced in [4], permits a natural characterization of the cellular space *C*.

*Definition 4.* A) The *quality function* defined with respect to $S_h^N$ is the map $\underline{b} : X \to \{-1, 0, 1\}^N$, with the value

$$\underline{b}(x) = [sgn(h_1(x)) \ldots sgn(h_N(x))], \ h_i \in S_h^N, \ \forall \ i \in I_N \qquad (4)$$

B) The value $\underline{b}(x)$, $x \in X$ is *consistent* if $sgn(h_i(x)) \neq 0$, $\forall \ h_i \in S_h^N$ and *inconsistent* else.□

The quality of the continuous state trajectory $x(.)$ can be defined as the composition of the functions $\underline{b}$ and $x(.)$, $(\underline{b} \circ x) : \mathbb{R} \to \{-1, 0, 1\}^N$, $(\underline{b} \circ x)(t) = \underline{b}(x(t))$, $\forall \ t \in \mathbb{R}$.

*Remark 1.* a) In *definition 3*, $Q = \text{card}(C) = \text{card}(P)$ and it will result that $Q \leq 2^N$, where *N* is the total number of hypersurfaces. *Q* is a measure of the complexity of the DES-plant automaton. b) The label-function *et* is built bijective. c) $\underline{b}(x)$ is consistent $\Leftrightarrow x \notin \text{Ker}(h_i)$, $\forall \ h_i \in S_h^N$. □

Some special properties of $S_h^N$, *Rel*, *C* and of the function $\underline{b}$ are listed below. Their proofs follow directly from the definitions above and are presented in the annex.

P1. $\forall \ x_1, x_2 \in DX$, it is true that $(x_1, x_2) \in Rel \Leftrightarrow \underline{b}(x_1) = \underline{b}(x_2)$. □

P2. $\forall \ x_1, x_2 \in DX$, it is true that $(x_1, x_2) \in Rel \Leftrightarrow$ there is a unique *N*-tuple of halfspaces, denoted $(H_1, \ldots, H_N)$, s.t. the following conditions are satisfied : 2.1) $H_i = H_i^- \in SH$ or $H_i = H_i^+ \in SH$, $\forall \ i \in I_N$ ; 2.2) $H_i \neq H_j$, $\forall \ i, j \in I_N$, $i \neq j$ ; 2.3) $H = \bigcap_{i=1}^{N} H_i \neq \emptyset$ ; 2.4) $x_1 \in H$ and $x_2 \in H$.□

P3. $\forall \ c_q \in C$, there is a unique *N*-tuple of halfspaces from the event topology *SH*, denoted $(H_1, \ldots, H_N)$, satisfying the conditions 2.1), 2.2) and 2.3) from P2, s.t.

$$c_q = \bigcap_{i=1}^{N} H_i \qquad (5)$$

P4. 4.1) Consider $c_q \in C$ an arbitrary cell, $I_N$ the index set of $S_h^N$ and $I_Q = \{1, \ldots, Q\}$. Then $\underline{b}(x)$ is constant, $\forall \ x \in c_q$ and $\underline{b}(x) = \boldsymbol{b}_q = [b_q^1 \ldots b_q^N] \in \{-1, 1\}^N$ (i.e. $\forall \ q \in I_Q$, $b_q^i \in \{-1, 1\}$, $\forall \ i \in I_N$). 4.2) Consider $c_q, c_p \in C$ two arbitrary cells and $\boldsymbol{b}_q = \underline{b}(x)$, with $x \in c_q$ and $\boldsymbol{b}_p = \underline{b}(x)$, with $x \in c_p$. Then $c_q \neq c_p \Leftrightarrow \boldsymbol{b}_q \neq \boldsymbol{b}_p$. □

The proof of P4 is trivial. P3 shows that a cell is an intersection of a unique combination of *N* open halfspaces determined by $S_h^N$. P4 shows that the quality function $\underline{b}$ can take only $Q = \text{card}(C) = \text{card}(P)$ distinct consistent values and that, based on *definitions 2, 3*A) and *4*A), there is a one to one correspondence between these consistent values and the cells from *C*. Define the set of all consistent quality values

$$B = \{\boldsymbol{b}_1, \ldots, \boldsymbol{b}_Q\} \qquad (6)$$

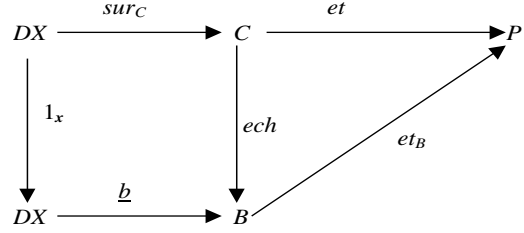

$$DX = X \setminus Fr, \quad Fr = \bigcup_{i=1}^{i_N} Ker(h_i)$$

Fig. 2 - The semantic content of the alphabet *P*.

with $\boldsymbol{b}_q = \underline{b}(\boldsymbol{x})$, $\boldsymbol{x} \in c_q \in C$, $\boldsymbol{x}$ arbitrary, $\forall q \in I_Q = \{1, \ldots, Q\}$. Two bijective maps can be defined,

$$ech : C \to B, \quad ech(c_q) = \boldsymbol{b}_q, \quad \forall q \in I_Q \tag{7}$$

and

$$et_B : B \to P, \quad et_B(\boldsymbol{b}_q) = p_q, \quad \forall q \in I_Q \tag{8}$$

with the property $et_B \circ ech = et$ (fig.2).

The next definition is a version of the concept of adjacency of two cells, introduced in [1].

*Definition 5*. Consider $c_a, c_b \in C$ and $\boldsymbol{b}_a = [b_a^1 \ldots b_a^N] = ech(c_a) \in B$, $\boldsymbol{b}_b = [b_b^1 \ldots b_b^N] = ech(c_b) \in B$. $c_a$, $c_b$ are *adjacent* if $\exists h_i \in S_h^N$ s.t. 5.1) $b_a^i b_b^i = -1$ and 5.2) $b_a^j b_b^j = 1$, $\forall j \in I_N, j \neq i$. $c_a$ and $c_b$ are called adjacent on $Ker(h_i)$. □

*Remark 2*. It's obvious that if $c_a$ and $c_b$ are adjacent on $Ker(h_i)$, then $h_i$ is unique, i.e. there is no other $h_q \in S_h^N$, $h_q \neq h_i$, s.t. $c_a$ and $c_b$ are adjacent also on $Ker(h_q)$. *Definition 5* states that the consistent quality values respectively associated to two cells which are adjacent on $Ker(h_i)$, differ only through their $i^{th}$ components. □

The next property is crucial for the observability of the DES-plant.

P5. If $c_a, c_b \in C$ are *adjacent* on $Ker(h_i)$, $h_i \in S_h^N$, then $c_b$ is the *unique* cell that is adjacent to $c_a$ on $Ker(h_i)$. □

## 3. PLANT-EVENTS AND PLANT-SYMBOLS

The general characterization of plant events, proposed in this paper is a variant of the definition introduced by Antsaklis [1]. Recall that $I_N = \{1, \ldots, N\}$ is the index set of $S_h^N$.

*Definition 6*. A) Consider a functional $h_i : X \to \mathbb{R}$, $h_i \in S_h^N$ and $\boldsymbol{x}(.)$ a continuous evolution of the plant model (1). The *plant-event* ($i+$) ( or ($i-$)) occurs at $t_e \in \mathbb{R}^+$ if the following conditions are satisfied : a) $h_i(x(t_e)) = 0$ ; b) $\exists\, d_1 > 0$ s.t. for all $\varepsilon$, $0 < \varepsilon < d_1$, $h_i(x(t_e+\varepsilon)) > 0$ (respectively $h_i(x(t_e+\varepsilon)) < 0$) ; c) $\exists\, d_2 > 0$ s.t. for all $\varepsilon$, $0 < \varepsilon < d_2$, $h_i(x(t_e-d_2)) < 0$ and $h_i(x(t_e-\varepsilon)) \leq 0$ (respectively $h_i(x(t_e-d_2)) > 0$ and $h_i(x(t_e-\varepsilon)) \geq 0$).

B) The set of all plant-events that can occur with respect to $S_h^N$ is $E_a = \{(i+), (i-) \mid i \in I_N\}$. The *alphabet of plant-symbols* is a set of $2N$ distinct symbols $Z = \{z_{i+}, z_{i-} \mid i \in I_N\}$, s.t. $\forall\, i \in I_N$, $z_{i+}$ (or $z_{i-}$) is associated to $(i+) \in E_a$ (or to $(i-) \in E_a$, respectively). □

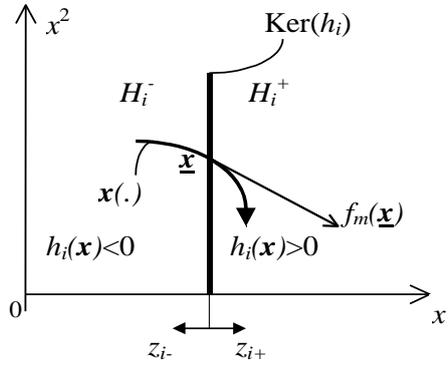

Fig. 3 - Example: the plant-event $(i+)$ occurs whenever the continuous state trajectory $x(.)$ crosses $\mathrm{Ker}(h_i)$ in the positive direction; $f_m(\underline{x}) = f(\underline{x},\, act(r_m))$, $r_m \in R$ and $x \in X = \mathbb{R}^2$.

According to *definition 6*, a plant event $(i+)$ ( or $(i-)$) occurs if (a) at *the time $t_e$ of the plant-event* the plants state lies on the triggering hypersuprface $\mathrm{Ker}(h_i)$; (b) *immediately after* the plant-event, the plant state lies on the positive side, included in $H_i^+$, for $(i+)$ (respectively on the negative side, included in $H_i^-$ for $(i-)$) of the triggering hypersurface $\mathrm{Ker}(h_i)$; (c) *prior* to reaching the triggering hypersuprface $\mathrm{Ker}(h_i)$, the plants state lied on the negative side, included in $H_i^-$, for $(i+)$ (respectively on the positive side, included in $H_i^+$, for $(i-)$) (fig.3).

The next assumption is generic and it means that a continuous state trajectory does not pass through the intersection point of two or more hypersurfaces from $S_h^N$.

*Hypothesis 1.* In a HCS, the plant-events do not occur simultaneously. □

Consider $Fr = \bigcup_{i=1}^{N} \mathrm{Ker}(h_i)$, $h_i \in S_h^N$, $\forall\, i \in I_N$ and define the set of all intersection points of the hypersurfaces of $S_h^N$, $Int = \{x \in Fr \mid \exists\, h_i, h_j \in S_h^N, h_i \neq h_j,$ s.t. $x = \mathrm{Ker}(h_i) \cap \mathrm{Ker}(h_j)\}$. In *hypothesis 1*, the mechanism of plant-symbols generation is described by the function $gev : (Fr \setminus Int) \to Z$

$$gev(x(t)) = \begin{cases} z_{i+},\, \text{if } (i-) \in E_a \text{ occures at } t \in \mathbb{R}^+ \\ z_{i+},\, \text{if } (i+) \in E_a \text{ occures at } t \in \mathbb{R}^+ \end{cases} \qquad (9)$$

A sequence of plant-events is denoted $w_e = e(1),e(2),…,e(k), …$ with $e(k) \in E_a$, $\forall\ k \in |N$, $k \neq 0$, and $t_e(k) \in |R^+$ the moment when $e(k)$ occurs, s.t. $t_e(k) < t_e(k+1)$, $\forall\ k \in |N, k \neq 0$ ; $t_e(0) = 0$ is the moment when the external event INIT occurs. $w_e$ produces a sequence of plant-symbols $w_z = z(1),z(2),…,z(k), …$ and $z(k) = gev(x(t_e(k))) \in Z$, $\forall\ k \in |N, k \neq 0$.

*Remark 3.* a) If *hypothesis 1* does not hold, then the generator in the interface can be designed in order to produce a sequence of plant symbols, s.t. the plant symbols are sent to the controller in the increasing order of their indexes $i \in I_N$.

b) There is a one to one correspondence between the plant-events set, the alphabet of plant-symbols and the event topology, i.e. $E_a \sim Z \sim SH$. □

## 4. THE DES-PLANT MODEL

### 4.1. DEFINITION AND SEMANTIC CONTENT

The continuous plant (1) coupled to the interface behaves like a DES representing all the possible logical evolutions of the quality $\underline{b}(x(.))$, from any initial state, $x(0) \in c_q$, $\forall\ c_q \in C$.

*Definition 7.* A *DES-plant model* is a nondeterministic automaton $G_p = \{P, R, f_p, Z, g_p\}$, where $P$ is the set of discrete states, $R$ is the input alphabet of control-symbols, $Z$ is the output alphabet of plant-symbols, $f_p : P \times R \to 2^P$ is the state transition function and $g_p : P \times P \to Z$ is the output function. The dynamic equations are :

$$p(k+1) \in f_p(p(k), r(k)), \quad g_p(p(k), p(k+1)) = z(k+1) \qquad (10)$$

where $p(k) \in P$, $z(k+1) \in Z$, $r(k) \in R$ and $p(k) \neq p(k+1)$, $\forall\ k \in |N$. □

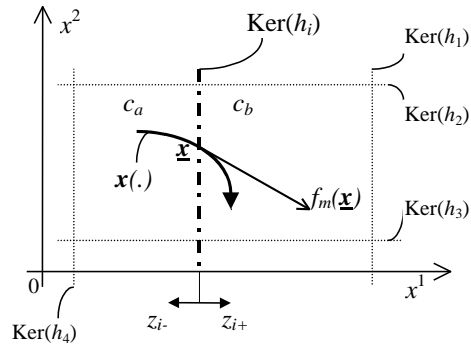

Fig. 4 – Example : the occurrence of the discrete state transition $p_a \to p_b$, observed as the plant-symbol $z_{i+}$, with $c_a = et^{-1}(p_a)$ and $c_b = et^{-1}(c_b)$ adjacent on $Ker(h_i)$; $f_m(\underline{x}) = f(\underline{x}, act(r_m))$, $r_m \in R$.

The definition of the output function $g_p$ reflects the *hypothesis 1* and differs from the one proposed in [1]. It's obvious that a discrete state transition, $p_a \rightarrow p_b$, $p_a, p_b \in P$ can occur in $G_p$ if and only if $\exists\ r_m \in R$ and $z_i \in Z$ s.t. $p(k) = p_a$, $p(k+1) = p_b$, $r(k) = r_m$, $z(k+1) = z_i$ satisfy (10).

Consider $sur_C$ the canonical surjection (*definition 3*.C) and fig.2), the plant (ODE) (1) and $k \in |N$ the present logical moment. Assume that at the moment $t_c(k) \in |R^+$, the continuous state is $x(t_c(k)) = x_0 \in c(k) \in C$, $G_p$ is in the discrete state $p(k) = et^{-1}(c(k))$ and receives the control-symbol $r(k) \in R$. By integrating the differential system (ODE) with the constant input $act(r(k)) \in U$, the continuous state at $t \geq t_c(k)$ is

$$F_k(\mathbf{x}_0, t) = \mathbf{x}_0 + \int_{t_c(k)}^{t} f(\mathbf{x}(\theta), act(r(k))) d\theta \tag{11}$$

The cell associated to $p(k)$ is $c(k) = sur_C(x_0)$. At the moment

$$\underline{t} = \min_{t}\ \{\ t > t_c(k)\ |\ et^{-1}(p(k)) \neq sur_C(F_k(x_0, t)),\ \text{with}\ F_k(x_0, t_c(k)) = x_0\ \} \tag{12}$$

the continuous state trajectory enters the cell $c(k+1) = sur_C(F_k(x_0, \underline{t})) \in C$. In *hypothesis 1*, $c(k+1)$ is *adjacent* with $c(k)$, on some hypersurface Ker($h_i$), $h_i \in S_h^N$ (fig.4). In conclusion, given $p(k) \in P$ and $r(k) \in R$, the discrete dynamics of $G_p$ from the logical moment $k$ to $k+1$ is described by the equations :

$$p(k+1) \in f_p(p(k), r(k)), \quad p(k+1) = et(sur_C(F_k(x_0, \underline{t}))),\ x_0 \in et^{-1}(p(k)) \tag{13.1}$$
$$g_p(p(k), p(k+1)) = z(k+1),\ z(k+1) = gev(F_k(x_0, t_e(k+1))),\ t_e(k+1) = \underline{t}\text{-}\varepsilon \tag{13.2}$$

So $G_p$ is in $p(k+1)$ if $\lim_{\varepsilon \rightarrow 0} x(t_e(k+1)+\varepsilon) \in et^{-1}(p(k+1))$, where $x(.)$ is the continuous evolution, $et^{-1}(p(k+1)) = c(k+1) \in C$ and $x(t_e(k+1)) \in \text{Ker}(h_i)$, with Ker($h_i$) the hypersurface that separates $c(k)$ from $c(k+1)$. $t_e(k+1) \in |R^+$ is the moment when the state transition $p(k) \rightarrow p(k+1)$ and the plant event $e(k+1) \in E_a$ occur, observed as the plant-symbol $z(k+1) \in Z$.

An evolution of $G_p$, $w_p = p(0), p(1), ..., p(k), ...$, $p(k) \in P$, $\forall\ k \in |N$, is controlled by a sequence of control-symbols $w_r = r(0), r(1), ..., r(k), ...$, $r(k) \in R$, $\forall\ k \in |N$ and is *observed* as a sequence of plant-symbols $w_z = z(1), z(2), ..., z(k+1), ...$, $z(k) \in Z$, $\forall\ k \in |N$.

### 4.2. THE OBSERVABILITY OF THE DES-PLANT AUTOMATON

<u>*Definition 8.*</u> The DES-plant automaton $G_p$ is *observable* if the current state can be determined uniquely from the previous state and plant-symbol, i.e. $\forall\ p_a, p_b, p_c \in P$ and $z_j \in Z$, if $z_j = g_p(p_a, p_b)$ and $z_j = g_p(p_a, p_c)$, then $p_b = p_c$.□

<u>*Proposition 1*.</u> The DES-plant automaton $G_p$ is *observable*.□

*Proposition 2*. If $G_p$ is observable, then for any initial state $p_0 \in P$ and sequence of admissible plant-symbols $w_z = z(1), z(2), ..., z(k+1), ..., z(k) \in Z, \forall k \in \mathbb{N}$, produced by the DES-plant automaton starting from $p_0$, there is a unique sequence of discrete states $w_p = p(0), p(1), ..., p(k), ...$, with $p(0) = p_0$, capable of producing $w_z$. □

## 5. EXAMPLE: THE OBSERVABILITY OF THE DES-PLANT AUTOMATON ASSOCIATED TO THE DOUBLE INTEGRATOR

The DES abstraction of the continuous double integrator is classic in the literature of hybrid systems [1], [2]. The extraction of the DES-plant automaton presented below was studied in [3]. The continuous double integrator, has the state equations

$$\dot{x}(t) = \begin{bmatrix} 0 & 1 \\ 0 & 0 \end{bmatrix} x(t) + \begin{bmatrix} 0 \\ 1 \end{bmatrix} u(t) \tag{14}$$

where $x = [x^1 \; x^2]^T \in \mathbb{R}^2$ is the state vector, $u \in \mathbb{R}$ is the control and $t \in \mathbb{R}^+$ is the continous time variable. The actuator can produce the discrete set of control values $U = \{-1, 0, 1\}$. The *alphabet of control-symbols* $R = \{r_1, r_2, r_3\}$ is associated to $U$ by the bijection $act : R \to U$, with the values: $act(r_1) = -1$, $act(r_2) = 0$, $act(r_3) = 1$. The corresponding state trajectories of the system (14) are defined by

$$x(t) = \begin{bmatrix} 1 & t \\ 0 & 1 \end{bmatrix} x_0 + \begin{bmatrix} 0.5t^2 \\ t \end{bmatrix} act(r_m), \quad r_m \in R, \quad x_0 = x(t_0) \in \mathbb{R}^2 \tag{15}$$

The equations of the integral curves are: (C1) $x^1(t) + (x^2(t))^2/2 = K$, for $act(r_1)$; (C2) $x^2(t) = K$, for $act(r_2)$; (C3) $x^1(t) - (x^2(t))^2/2 = K$, for $act(r_3)$. $K \in \mathbb{R}$ is a constant depending on the initial state $x_0$.

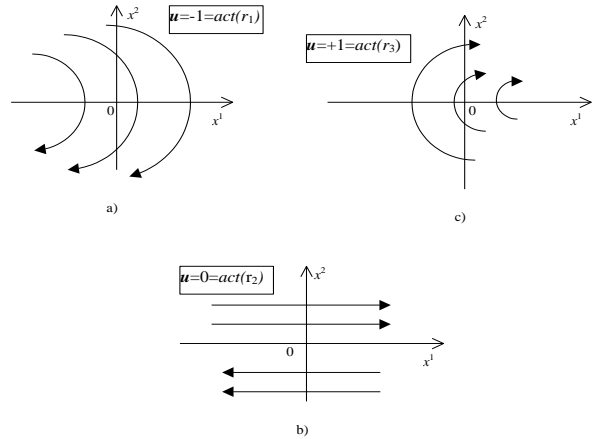

Fig. 5 - The state trajectories (15) of the double integrator (14).

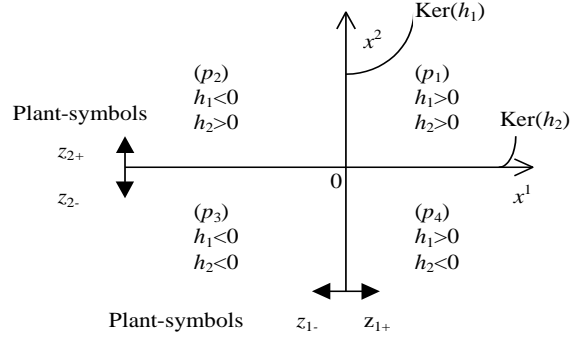

Fig. 6 - The state-space partition determined by $h_1(\boldsymbol{x}) = x^1$ and $h_2(\boldsymbol{x}) = x^2$, the alphabet $P$ of discrete states and the plant-symbols for the double integrator.

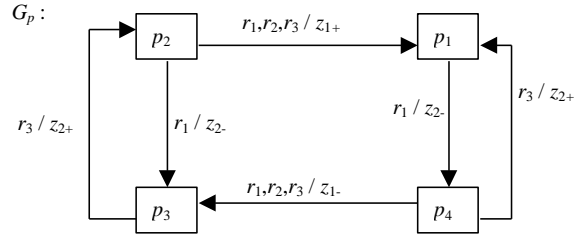

Fig. 7 - The DES plant model of the double integrator with the state space partition $S_h^2$.

*Table 1.*

| The alphabet $P$ | The cellular space $DX/_{Rel} = C$ | The set $B$ |
|---|---|---|
| $p_1$ | $c_1 = H_1^+ \cap H_2^+$ | $\boldsymbol{b}_1 = [1 \quad 1]$ |
| $p_2$ | $c_2 = H_1^- \cap H_2^+$ | $\boldsymbol{b}_2 = [-1 \quad 1]$ |
| $p_3$ | $c_3 = H_1^- \cap H_2^-$ | $\boldsymbol{b}_3 = [-1 \quad -1]$ |
| $p_4$ | $c_4 = H_1^+ \cap H_2^-$ | $\boldsymbol{b}_4 = [1 \quad -1]$ |

Consider the set of two functionals $S_h^2 = \{h_i : \mathbb{R}^2 \to \mathbb{R} \mid h_1(\boldsymbol{x}) = x^1, h_2(\boldsymbol{x}) = x^2\}$. Each functional $h_i \in S_h^2$ *separates* $\mathbb{R}^2$ into two open halfspaces $H_i^- = \{\boldsymbol{x} \in \mathbb{R}^2 \mid h_i(\boldsymbol{x}) < 0\}$, and $H_i^+ = \{\boldsymbol{x} \in \mathbb{R}^2 \mid h_i(\boldsymbol{x}) > 0\}$. The hypersurfaces $\mathrm{Ker}(h_i) = \{\boldsymbol{x} \in \mathbb{R}^2 \mid h_i(\boldsymbol{x}) = 0\}$, $i=1,2$ (fig.6) are nonsingular. *The quality function* is $\underline{b} : \mathbb{R}^2 \to \{-1, 0, 1\}^2$, $\underline{b}(\boldsymbol{x}) = [sgn(h_1(\boldsymbol{x})) \; sgn(h_2(\boldsymbol{x}))]$. Define $DX = \mathbb{R}^2 \setminus Fr$ with $Fr = \bigcup_{i=1}^{2} \mathrm{Ker}(h_i)$. The cellular space $DX/_{Rel} = C = \{c_1, c_2, c_3, c_4\}$, the set $B$ of all *consistent quality values* and the alphabet $P$ of discrete states are described in *Table 1*. For example: $p_1$ is the symbol of the open cell $c_1 = H_1^+ \cap H_2^+$ and of $\boldsymbol{b}_1 = [1 \quad 1] = \underline{b}(\boldsymbol{x})$, $\forall \boldsymbol{x} \in c_1$.

From the phase-portraits in fig.5 and from the state space partition (fig.6) the DES-plant automaton in fig.7 is deduced. $G_p$ is observable. For example, consider $p(0) = p_2$ and the sequence of plant-symbols $w_z = z(1),...,z(k), ... = z_{1+}z_{2-}z_{1-}z_{2+}z_{1+} ... = (z_{1+}z_{2-}z_{1-}z_{2+})^*$. The pair $(p_2, w_z)$ corresponds uniquely to the admissible discrete evolution $w_p = p(0),p(1),...,p(k), ... = p_2p_1p_4p_3p_2 ... = (p_2p_1p_4p_3)^*$. Concretely, the transition $p_2 \to p_1$ is observed as $z_{1+}$ and, given $p_2$ and $z_{1+}$, the next state $p_1$ results uniquely; the transition $p_1 \to p_4$ is observed as $z_{2-}$ and so on.

*Remark 4.* The DES-plant automaton extracted in [1] is also observable, but it uses four distinct functionals to define the state space partition and another mechanism for plant-symbols generation. □

## CONCLUDING REMARKS

Based on a version of the Antsaklis framework for HCS, the properties of the cellular space have been investigated and a mechanism for plant-symbols generation has been proposed. Property P5 states that if two cells are adjacent on a hypersurface, then there is no other cell being adjacent with one of the two cells on the same hypersurface. Assuming that the plant-events do not occur simultaneously in a HCS, it was shown that the previous property ensures the observability of the DES-plant automaton. The observability property is crucial for the proper decision taken by the DES-controller.

## ANNEX

P1. Proof. Consider $x_1, x_2 \in DX$, s.t. $(x_1, x_2) \in Rel$. Then $(x_1, x_2) \in Rel \Leftrightarrow h_i(x_1)h_i(x_2) > 0, \forall h_i \in S_h^N \Leftrightarrow sgn(h_i(x_1)) = sgn(h_i(x_2)) \neq 0, \forall h_i \in S_h^N \Leftrightarrow \underline{b}(x_1) = \underline{b}(x_2)$ and, according to *definition 4*.B), the quality value is consistent. □

**P2. Proof.** "$\Rightarrow$" Assume that $(x_1, x_2) \in Rel$. According to P1, it results that $\underline{b}(x_1) = \underline{b}(x_2) \in \{-1, 1\}^N$. Denote $\boldsymbol{b}_q = [b_q^1 \ldots b_q^N] = \underline{b}(x_1) = \underline{b}(x_2)$. Because $\boldsymbol{b}_q$ is a consistent value, it results that $\forall\, i \in I_N$ either $b_q^i = -1$ or $b_q^i = 1$, i.e. either $x_1, x_2 \in H_i^- \in SH$, or $x_1, x_2 \in H_i^+ \in SH$. Without loss of generality, assume that $b_q^i = 1$, $\forall\, i \in I_N$, which means that $x_1, x_2 \in H_i^+$, $\forall\, i \in I_N$. It results that $\exists$ an $N$-tuple of distinct halfspaces from $SH$, $(H_1^+, \ldots, H_N^+)$, s.t. $H = \bigcap_{i=1}^{N} H_i \neq \varnothing$ and $x_1 \in H$ and $x_2 \in H$. Moreover, the $N$-tuple of halfspaces from $SH$ is unique, because the value $\boldsymbol{b}_q$ is unique.

"$\Leftarrow$" Consider $x_1, x_2 \in DX$ and $(H_1^+, \ldots, H_N^+)$ an $N$-tuple of halfspaces from $SH$ s.t. the conditions 2.2)-2.4) are fulfilled. $x_1 \in H$ implies that $x_1 \in H_i^+$, $\forall\, i \in I_N$. It results that $h_i(x_1) > 0$, $\forall\, h_i \in S_h^N$ and $\underline{b}(x_1) = [sgn(h_1(x_1)) \ldots sgn(h_N(x_1))] = [1 \ldots 1] = \boldsymbol{b}_q$. Similarly, $x_2 \in H$ implies that $\underline{b}(x_2) = [1 \ldots 1] = \boldsymbol{b}_q$. Hence $\underline{b}(x_1) = \underline{b}(x_2) = \boldsymbol{b}_q$ and $\boldsymbol{b}_q$ is consistent so, according to P1, it results that $(x_1, x_2) \in Rel$. The proof is similar for any N-tuple of halfspaces that satisfies 2.1), 2.2) and 2.3). $\square$

**P3. Proof.** Consider $c_q \in C$ an arbitrary cell and $x_1, x_2 \in c_q$, $x_1, x_2$ arbitrary. According to *definition 3*.A), $(x_1, x_2) \in Rel$. According to P2, $(x_1, x_2) \in Rel \Leftrightarrow \exists \cup (H_1, \ldots, H_N)$ an N-tuple of distinct halfspaces from $SH$ s.t. the conditions 2.1)-2.4) in P2 are satisfied. Because $c_q$ is a class of equivalence it results that $(x_1, x_2) \in Rel$, $\forall\, x_1, x_2 \in c_q$, so $c_q = \bigcap_{i=1}^{N} H_i$. $\square$

**P5. Proof.** Consider $\boldsymbol{b}_a = [b_a^1 \ldots b_a^N] = ech(c_a) \in B$ and $\boldsymbol{b}_b = [b_b^1 \ldots b_b^N] = ech(c_b) \in B$. Assume that $\exists\, c_c \in C$, $c_c \neq c_b$ s.t. $c_c$ is adjacent to $c_a$ on $Ker(h_i)$ and consider $\boldsymbol{b}_c = ech(c_c) = [b_c^1 \ldots b_c^N] \in B$. From the adjacency of $c_a$ and $c_b$ it results : (1) $b_a^i b_b^i = -1$ and $b_a^j b_b^j = 1$, $\forall\, j \in I_N, j \neq i$. From the adjacency of $c_a$ and $c_c$ it results : (2) $b_a^i b_c^i = -1$ and $b_a^j b_c^j = 1$, $\forall\, j \in I_N, j \neq i$. From (1) and (2) and from the fact that $b_a, b_b, b_c \in \{-1, 1\}^N$, it results that $\boldsymbol{b}_c = \boldsymbol{b}_b$. Because *ech* is bijective, it follows that $c_c = c_b$, which contradicts the initial assumption. Hence, $c_b$ is unique.$\square$

*Proposition 1.* Proof. Assume that $\exists\, p_a, p_b, p_c \in P$ and $z_j \in Z$ s.t. $z_j = g_p(p_a, p_b)$, $z_j = g_p(p_a, p_c)$ and $p_b \neq p_c$. Without loss of generality, assume that $z_j = z_{i+}$, $i \in I_N$. If $z_{i+} = g_p(p_a, p_b)$, then $\exists\, r_m \in R$ s.t. $p_b \in f_p(p_a, r_m)$, so there is a possible state transition $p_a \rightarrow p_b$, and, in *hypothesis 1*, this means that the cells $c_a = et^{-1}(p_a)$ and $c_b = et^{-1}(p_b)$ are *adjacent* on the hypersurface $Ker(h_i)$, $h_i \in S_h^N$. From $z_{i+} = g_p(p_a, p_c)$, in *hypothesis 1*, it results similarly that $c_a = et^{-1}(p_a)$ and $c_c = et^{-1}(p_c)$ are also adjacent on $Ker(h_i)$. From P5 it results that $c_b = c_c$, hence $p_b = p_c$ and the assumption "$p_b \neq p_c$" is contradicted. In consequence, $G_p$ is observable.$\square$

*Proposition 2.* Proof. The definition of the observability can be applied iteratively to prove that each state $p(k+1)$ of the sequence $w_p$ is determined, uniquely by the previous state $p(k)$ and current plant-symbol $z(k+1)$. $\square$